\documentclass[12pt]{article}

\usepackage[english]{babel}
\usepackage[usenames]{color}
\usepackage[cp1250]{inputenc}
\usepackage{amsfonts}
\usepackage{amsthm}
\usepackage{graphicx}
\usepackage{epsfig}

\theoremstyle{plain}
\newtheorem{df}{Definition}

\newtheorem{sps}[df]{Observation}

\begin{document}
\newcommand{\bea}{\begin{eqnarray}}
\newcommand{\eea}{\end{eqnarray}}
\newcommand{\be}{\begin{equation}}
\newcommand{\ee}{\end{equation}}
\newcommand{\beas}{\begin{eqnarray*}}
\newcommand{\eeas}{\end{eqnarray*}}
\newcommand{\bs}{\backslash}
\newcommand{\bc}{\begin{center}}
\newcommand{\ec}{\end{center}
}

\title{Analysis of the convex hull of the attractor of an IFS}

\author{Jarek Duda}
\date{\it \footnotesize Jagiellonian University, Reymonta 4, 30-059 Kraków, Poland, \\
\textit{email:} dudaj@interia.pl}

\maketitle

\begin{abstract}
In this paper we will introduce the methodology of analysis of the
convex hull of the attractors of iterated functional systems (IFS) -
compact fixed sets of self-similarity mapping: \be
K=\bigcup_{i=1..n}A_i(K)+t_i\ee where $A_i$ are some contracting,
linear mappings. The method is based on a function which for a
direction, gives width in that direction. We can write the self
similarity equation in terms of this function, solve and analyze
them. Using this function we can quickly check if the distance from
$K$ of a given $x$ is smaller than a given distance or even compute
analytically convex hull area and the length of its boundary.
\end{abstract}

\section{Introduction}

Fix $X=\mathbb{R}^m$ with $\|x\|:=\sqrt{x^Tx}$ norm,
$A_i:X\rightarrow X (i=1..n)$ are contracting matrices
$$c_i:=\|A_i\|:=\sup_{x:\|x\|=1}\|A_ix\|<1$$ and
$t_i\in X$ for $i=1..n$ are translations.\\

The finite set of contracting mappings $\{x\rightarrow
A_ix+t_i\}_{i=1..n}$ is called iterated functional system. In
\cite{hut} it is proved that it has an attractor: a unique compact
nonempty set K, which is fixed point of $I$:\be \label{ss1}
I(K):=\bigcup_{i=1..n} \left( A_i(K)+t_i \right)\ee
where $A(K):=\{Ax:x\in K\}$.\\

 A precise analysis of this attractor is usually very
difficult. It's common problem in e.g. computer graphics, image
compression to find a good approximations from above of this set.
There are known approaches to bound this set using spheres:
Rice\cite{rice}, Hart and DeFanti\cite{har}, Sharp and
Edalat\cite{sharp} or
boxes \cite{chu}.\\
In this paper will be introduced the methodology to bound it with
convex
set too, but this time - optimally: we will show how to construct its convex hull
and how to use to get quickly as precise approximation of the attractor as needed.\\
We will reduce this problem to the solution of $n-1$ dimensional
functional equation for the width function, which can be easily
approximated numerically and in some cases even found
analytically.\\

In Section 2 we will define the width function - it's a function
which for a given direction, gives a position of orthogonal
hyperplane bounding the set in this direction. This function
completely defines the convex hull of
closed set.\\
To describe a convex set, it's better to use the radius function for
given point, which gives for a given direction, length of segment
from that point in that direction. We will show how to get the radius function from the width function.\\
In Section 3 we will show how to change the self similarity equation
for $K$ to functional equation of its width function and that we can
solve it numerically by a simple iteration. \\
In Section 4 we will show how to use the width function in practice
- for example to decide if the distance from the point $x$ to
$K$ is smaller than a given number. \\
In Section 5 we will show explicit formula for the width function
for IFS with all
$A_i$ equal.\\
Then we will concentrate on simple 2-dimensional IFS. We will show
that it has a point of symmetry and using the width function - that
its convex hull is built of triangles. So we can compute the area of
the convex hull and the length of its boundary. \\
Using the isodiametric inequality we can get an interesting
trigonometric inequality:
$$\forall_{1<r\in \mathbb{R},\ \phi\in[0,2\pi]}\quad
\sum_{j>0}|\sin(j\phi)|r^{-j}\leq\frac{1}{\pi}\frac{r+1}{r-1}.$$

\section{Width function}

\begin{df}
  \emph{For a bounded, nonempty set $K\subset X$ we define }width function $h_x:D\to
  \mathbb{R}$ around $x\in X$:
  \beas h_{K,x}(d):=\inf\{h:K\subset x+H(d,h)\} \eeas
\end{df}
\begin{tabular}{lll}
where & $D:= S_{m-1}=\{x\in \mathbb{R}^m :\ x^Tx=1\}$ &-
directions,\\
&$H(d,h):=\{x\in \mathbb{R}^m:x^Td\leq h\}$ &- halfplane.
\end{tabular}\\

$K$ will be usually fixed, so we will write $h_x\equiv h_{K,x}$\\
Of course $h_x(d)=h_y(d)+(y-x)^Td$\\
If $x_0$ is in the convex hull of $K$, $h_{x_0}$ is nonnegative.\\
Obviously $h$ is bounded and it is easy to show \cite{bus} that is
continuous, and that it completely describes the convex hull of
compact set: \be L\equiv\textrm{conv}(K)=\bigcap_{d\in
D}H(d,h_{K,0}(d)).\ee

Fix a nonempty, compact, convex, set $L\subset X$ and $x_0\in L$.
\begin{df}
  \emph{For $L,\ x_0$ as above we define }radius function $r:D\to
  \mathbb{R}^+\cup\{0\}$ around $x_0$:
  \be (r_{L,x_0}(d)\equiv)\ r(d):=\sup\{r:x_0+rd\in L\} \ee
\end{df}

Now we can analyze the relation between the radius and the width
functions of $L$
around the same point, say $x_0=0$ ($h\equiv h_0$).\\
For any $s\geq 0$ \beas r(d)\leq s \Leftrightarrow \forall_{e\in
D:d^Te>0}\ sd\in H(e,h(e))\Leftrightarrow\forall_{e\in D:d^Te>0}\
sd^Te\leq h(e)\Leftrightarrow\forall_{e\in D:d^Te>0}\ s\leq
\frac{h(e)}{d^Te}\eeas So \be \label{e4} r(d)=\inf_{e\in
D:e^Td>0}\frac{h(e)}{d^Te}\ee

If $e$ fulfils infinium for $d$ we say that $e$ \emph{supports}
$d$.\\

We would like to use differential methods - we have to expand $h$
for a moment in a neighborhood of the sphere. We would do it in the
simplest way - take $h(\hat{e})$ to ensure that $$e^T\nabla
h(\hat{e})=0$$ where $\hat{x}:=x/\| x\|, \quad \nabla\equiv
(\partial_{e_i})_{i=1..m}$.

Assume now that $h(\hat{e})$ is differentiable for some $e\in D$.\\
So the necessity condition for $e\in D$ to support some $d\in D$
from (\ref{e4}) is
$$\nabla\frac{h(\hat{e})}{d^Te}\in \mathbb{R}e\Leftrightarrow \exists_{\lambda\in\mathbb{R}}\ (d^Te)\nabla h(\hat{e})-h(e)d=\lambda e$$
multiplying by $e^T$, we get $\lambda=-h(e)(e^Td)$

 $$ d=(d^Te)\left(\frac{\nabla
 h(\hat{e})}{h(e)}+e\right)$$
Now we take $d^Te\geq 0$ to fulfil $\|d\|=1$, (\ref{e4}) becomes
 \be \label{nec} r\left(\widehat{\frac{\nabla
 h(\hat{e})}{h(e)}+e}\right)=h(e)\left\|\frac{\nabla
 h(\hat{e})}{h(e)}+e\right\|. \ee

\section{Self similarity relation}
For any nonsingular matrix $A:X\to X$, translation $t\in X$, $x\in
X$, $d\in D$, $a\in\mathbb{R}^+$ \be \label{s31} Ax+t\in H(d,a)
\Leftrightarrow
 x^TA^Td\leq a-t^Td\Leftrightarrow x\in H\left(\widehat{A^Td},\frac{a-t^Td}{\|
A^Td\|}\right) \ee
Now look at the self-similarity equation (1):
\beas & h_0(d)\leq a\Leftrightarrow K\subset H(d,a)\Leftrightarrow
\forall_i\ A_iK+t_i\subset H(d,a)
\Leftrightarrow^{(\ref{s31})}\forall_i\ K\subset
H\left(\widehat{A_i^Td},\frac{a-t_i^Td}{\| A_i^Td\|}\right)\Leftrightarrow\\
&\Leftrightarrow \forall_i\ h_0(\widehat{A_i^Td})\leq
\frac{a-t_i^Td}{\| A_i^Td\|}\Leftrightarrow
 \max_i\left(
\|A_i^Td\|h_0(\widehat{A_i^Td})+t_i^Td\right)\leq a \eeas It's true
for any $a\in\mathbb{R}^+$, so we've get functional equation for the
width function:

\begin{sps}
\be \label{ss2} h_0(d)=\max_i
\|A_i^Td\|h_0(\widehat{A^Td})+t_i^Td.\ee
\end{sps}

In some cases we can solve this equation analytically, but numerical
approximation should usually be enough.

Consider $\mathbf{C}(D,\mathbb{R})$ - the space of continuous
functions $D\to\mathbb{R}$ with supremum norm: \\$\|f\|=\sup_{d\in
D}f(d)$.\\
Define $I:\mathbf{C}(D,\mathbb{R})\to \mathbf{C}(D,\mathbb{R})$ \be
\label{it} I_h(f)(d):=\max_i\left(
\|A_i^Td\|f(\widehat{A_i^Td})+t_i^Td\right)\ee Now for any $f,g\in
\mathbf{C}(D,\mathbb{R})$: \setlength\arraycolsep{2pt} \beas
\|I_h(f)-I_h(g)\|&\leq&\sup_{d\in
D}\sup_{r\in[-R,R]}\left|\max_i\left[\|A_i^Td\|(g(\widehat{A_i^Td})+r)+t_i^Td\right]-
\max_j\left[\|A_j^Td\|g(\widehat{A_j^Td})+t_j^Td\right]\right|\leq \\
&\leq& \sup_{d\in
D}\sup_{r\in[-R,R]}\left|\max_i\left(\|A_i^Td\|r\right)\right|
 = R\max_i\|A_i\|=\|f-g\|c\eeas
\begin{tabular}{ll}where & $R:=\|f-g\|$,\\ & $c:=\max_i\|A_i\|=\max_i c_i<1.$ \end{tabular}

So because $A_i$ are contracting, $I_h$ is contracting with $c$
coefficient.\\

 Using Banach contraction theorem, we get the unique fixed point
of iteration (\ref{it}) - the width function of our attractor.\\
So to approximate the width function numerically, we can start from
a constant function (the width function of a ball) and iterate
(\ref{it}).
\section{Approximation of attractor}
In this section it will be shown how to use found width function to
approximate $K$ as precise as needed.\\
To check if $x\in L$ (convex hull of $K$) we should check if
$$x\in L\leftrightarrow\forall_{e\in D} (x-x_0)^Te\leq h_{x_0}(e)
\leftrightarrow \|x-x_0\|\leq r_{x_0}(\widehat{x-x_0})$$ for some
fixed $x_0\in L$.\\
We can check it immediately having the radius function, but finding
this function from the width function(\ref{nec}), requires some
smoothness - can be generally
difficult, especially for numerical approximated functions.\\
We will see that we won't loose much of precision, if instead of
checking all directions for the width function, we will check only
one: $\widehat{x-x_0}$
$$L_{x_0}:=\{x:(x-x_0)^T(\widehat{x-x_0})\leq h_{x_0}(\widehat{x-x_0})\}=
\{x:x^T(\widehat{x-x_0})\leq h_0(\widehat{x-x_0})\}$$
\begin{figure}[h]
    \centering
        \includegraphics[width=15cm]{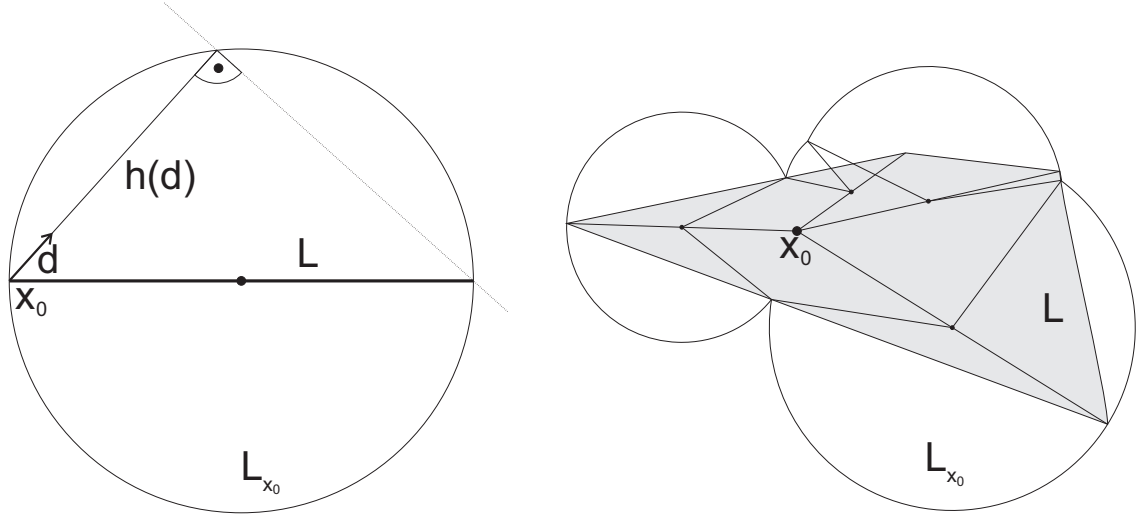}
        \caption{$L_{x_0}$ for a segment (left) and some polygon.}
\end{figure}
From fig. 1 we see that $L_{x_0}$ for a segment, where $x_0$ is one
of its ending points, is a ball with the center in the middle of
this segment.\\
We can threat convex set($L$) as the sum of such segments for all
directions, so $L_{x_0}$ is the sum of all such balls (in fact we
can restrict to the supporting points of $K$).\\
So we have rough approximations:
$$L\subset L_{x_0}\subset B(x_0,R_{x_0})$$
$$\rho(L_{x_0},L)\leq R_{x_0}/2$$
\begin{tabular}{ll}where & $\rho(A,B):= \sup_{x\in A}\inf_{y\in B}\|x-y\|$,\\ &
$R_{x_0}:=\sup_{d\in D}h_{x_0}(d)\ (=\inf_{r>0}\{L\subset
B(x_0,r)\})$. \end{tabular}\\
If $h_{x_0}$ isn't constant ($L$ is a ball) this approximation is
better than with a ball.\\

We could go closer to $L$ this way, by checking more points, but we
can use this additional checkings to came closer to the attractor -
check if $x\in I^k(L_{x_0})$:
$$x\in I^{l+1}(L_{x_0})\Leftrightarrow \exists_i\ x\in A_i(I^l(L_{x_0}))+t_i\Leftrightarrow
\exists_i\ A_i^{-1}(x-t_i)\in I^l(L_{x_0})$$

If we iterate this equivalence $k$ times, we get an algorithm
checking if $x\in I^k(L_{x_0})$:\\

\verb"function near("$x$\verb",k)"\\
\verb"     { int i=0;"\\
\verb"       if ("$x^T(\widehat{x-x_0})>h_0(\widehat{x-x_0})$\verb") return false else"\\
\verb"        if(k=0) return true else"\\
\verb"         { do i++ until ((i>n) or ("$x\in \textrm{Im}(A_i)$\verb" and near("$A_i^{-1}(x-t_i)$\verb",k-1)));"\\
\verb"           if(i>n) return false else return true"\\
\verb"         }"\\
\verb"     }"\\

The condition with the image of $A_i$ is required for singular
matrices and numerically cannot be fulfilled - we can just omit this matrices.\\

To evaluate the results of this algorithm, define
$$C_k=\rho(I^k(L_{x_0}),K)$$
$C_0\leq R_{x_0}/\sqrt{2}$ is some constant,
which should be approximated from above eg numerically.\\

Let's make the iteration now: \be\label{s32}\inf_{y\in
K}\|x-y\|=\inf_i\inf_{y\in
K}\|x-t_i-A_i(y)\|\leq\inf_i\|A_i\|\inf_{y\in
K}\|A_i^{-1}(x-t_i)-y\|\ee
So $C_k\leq C_{k-1}\max_i\|A_i\|= C_{k-1}\max_ic_i\leq C_0c^k$.\\

Now using (\ref{s32}) we can alter this algorithm: if
\verb"near1("$x,l$\verb")" will return \verb"true", we are sure that
$x$ is nearer
to $K$ than a given distance $l$:\\

\verb"function near1("$x,l$\verb")"\\
\verb"     { int i=0;"\\
\verb"       if ("$x^T(\widehat{x-x_0})>h_0(\widehat{x-x_0})$\verb") return false else"\\
\verb"        if("$l\geq C_0$\verb") return true else"\\
\verb"         { do i++ until ((i>n) or near1("$A_i^{-1}(x-t_i),l/c_i$\verb"));"\\
\verb"           if(i>n) return false else return true"\\
\verb"         }"\\
\verb"     }"\\

Analogically we can construct for example an algorithm to answer
questions like if $I^k(L)\bigcap (I^k(L)+t)=\emptyset$.

\section{Analytically solvable examples}
In this section we will show that in some cases we can compute width
function analytically.

\begin{sps}In the case $A_1=A_2=...=A_n(=A)$ the width function realize
  \be \label{ss3} h_0(d)=\|Ad\|h_0(\widehat{Ad})+h^*(d)\ee
  where $h^*(d)=\max_i{t_i^Td}=h_{\{t_1,t_2,...,t_n\},0}(d)$\\
  Solution to this equation is:
  \be \label{sol} h_0(d)=\sum_{i=0}^\infty\|A^id\|h^*(\widehat{A^id})\ee
\end{sps}
\textbf{Proof:} (\ref{ss3}) is obvious from (\ref{ss2})\\
By induction over $k$, using (\ref{ss3}): \beas
h_0(d)=\|A^kd\|h_0(\widehat{A^kD})+\sum_{i=0}^{k-1}\|A^id\|h^*(\widehat{A^id})\eeas
$A$ is contracting, $h$ is bounded - the first term tends to 0. $\qquad\square$\\

Now we will show an example how to get complete set of
information about the convex hull $L$ of $K$ using the width function.\\

Fix $X=\mathbb{R}^2$. We can threat $X=\mathbb{C}$ as the complex
plane,
multiplication by a complex number corresponds to a rotation and a scaling.\\
We can identify the directions space $D$ with angles:
$\alpha\equiv(\cos\alpha,\sin\alpha)$.\\

Fix $re^{\mathbf{i}\phi}=z\in\mathbb{C}:r=|z|>1,\ 2\leq n\in \mathbb{N}$\\
In the rest of this section we will analyze attractor, which can
represent fractional part in a complex base system\cite{me}: \be
\label{self} Kz=\sum_{i=0..n-1}\left(K+i\right)\ee

Before investigating the width function, we observe that $K$ has the
symmetry point:
\begin{sps}
$x_0=\frac{1}{2}\frac{n-1}{z-1}$ is the center of symmetry of $K$.
\end{sps}
\textbf{Proof:} For any point $x_0$, from (\ref{self}) $2x_0-K$
satisfies:
$$z(2x_0-K)=2x_0z-\sum_{i=0..n-1}(K+i)=\sum_{k=n-1-i=0..n-1}2x_0z-n+1-K+k$$
Hence if $\ 2x_0z-n+1=2x_0$ then $2x_0=\frac{n-1}{z-1}$ and we have:
$$z(2x_0-K)=\sum_{i=0..n-1}(2x_0-K)+i$$
Because of the uniqueness: $2x_0-K=K$, so $\ x_0-K=K-x_0$.
$\quad\quad\quad\Box$\\

We have $h_{x_0}(\alpha)\equiv
h(\alpha)=h(-\alpha)=(h_0(\alpha)+h_0(-\alpha))/2$\\
Now $h^*(\alpha)+h^*(-\alpha)=(n-1)r^{-1}|\cos(\alpha+\phi)|$,
(\ref{sol}) gives \be
2h(\alpha)=(n-1)\sum_{j>0}r^{-j}|\cos(\alpha+j\phi)|\ee

Now assume that $\phi=\pi \frac{l}{k}\quad (l,k\in\mathbb{N})$
\be\label{rat}2h(\alpha)=(n-1)\sum_{j=1..k}|\cos(\alpha+j\phi)|\sum_{i\geq
0}r^{-j-ki}=\frac{n-1}{1-r^{-k}}\sum_{j=1..k}r^{-j}|\cos(\alpha+j\phi)|\ee

In this case, our width function is differentiable everywhere except
finite number of points. We could use (\ref{nec}) ($\nabla
h(\widehat{\alpha})\equiv (-\sin\alpha,\cos\alpha)h'(\alpha)$) for
differentiable points of $h$ around indifferentiable, but instead we
will show that this points corresponds to straight lines on the
boundary of $L$ (convex hull of $K$). In
fact, it will occur that $L$ is a polygon.\\

\begin{figure}[h]
    \centering
        \includegraphics[width=18cm]{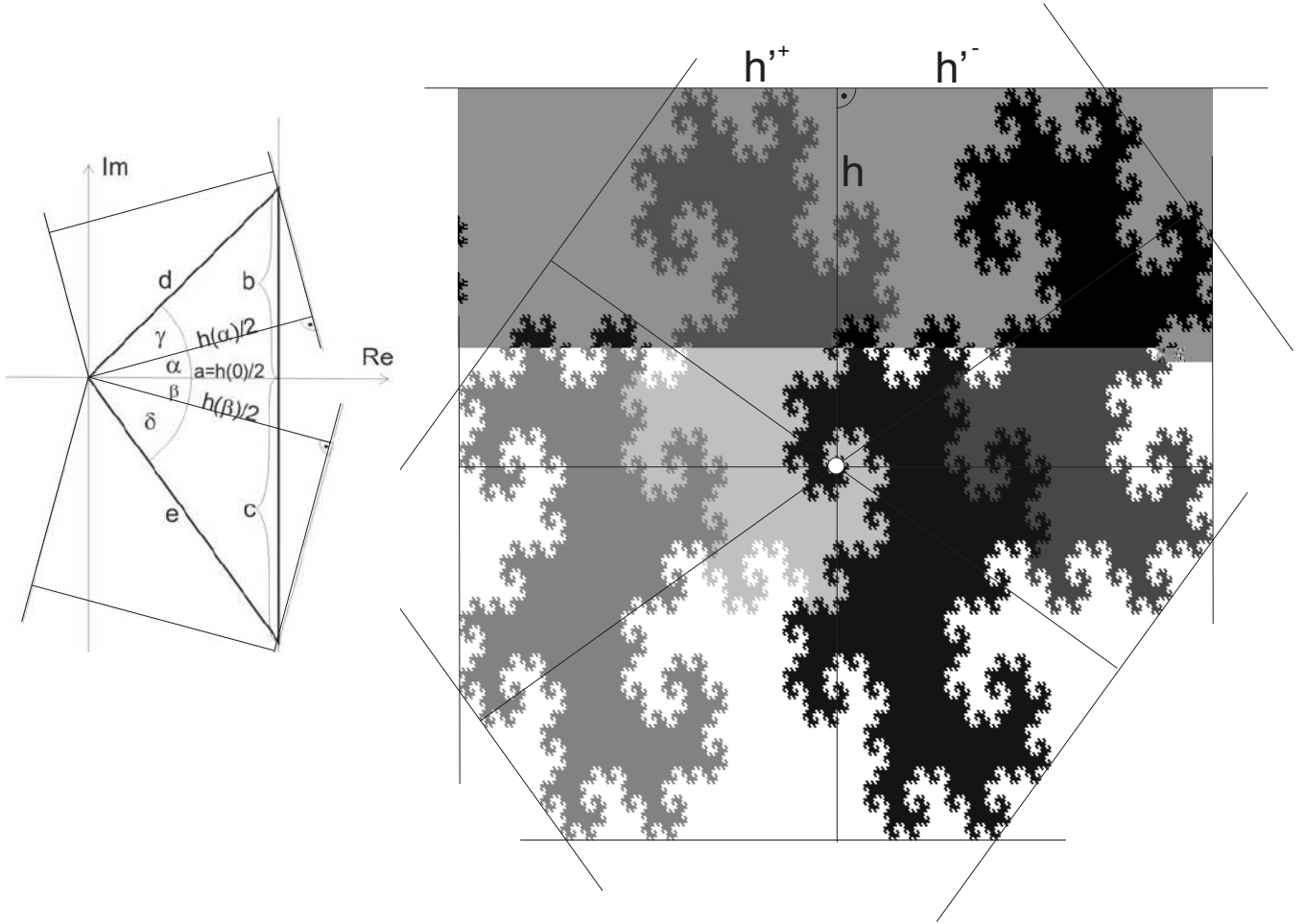}
        \caption{Analysis of the behavior of width function around a
        segment on the boundary of a convex set(left) and
        construction of the convex hull in case: $z=i+1,\ n=2$ (right).}
\end{figure}

Look at the triangle on  fig. 2 - width function will have
indifferentiable minimum in 0 - in a neighborhood of 0:
$$h(\alpha)=\bigg\{\begin{array}{ll}
               d\cos(\gamma)&0\geq\alpha<\epsilon \\
               d\cos(\delta)&0>\alpha>-\epsilon
             \end{array}
$$
\be\label{wid}h'(0^+)=d\sin(\alpha+\gamma)=b\quad,\quad
h'(0^-)=-e\sin(\beta+\delta)=-c\ee So we see:
\begin{enumerate}
  \item If a vertex in some direction $\gamma$ is supported by some range of
  directions - in that range $h(\alpha)=d\cos(\gamma-\alpha)$.
  Thanks to the uniqueness of $h$, implication above is equivalence.
  \item In a direction in which the edge of the convex hull contains
  straight segment, iff the width function have local minimum and is
  indifferentiable.\\
  We can construct this segment using (\ref{wid}).
\end{enumerate}
Now look at (\ref{rat}) - there is a finite number of
indifferentiability points and between them function is in the form
of $h(\alpha)=\sum_{i=0..k-1}a_i\cos(\alpha+b_i)$. Expanding cosinus
of sum, we see that it can be written in the form
$h(\alpha)=a\cos(\alpha+b)$ for some $a,b\in \mathbb{R}$. So it
supports on one point - common point of succeeding triangles -  our convex hull is polygon.\\

We will construct it now - call $T(a,b,c)$ - triangle as in fig. 2\\
This triangles should be rotated by angles of indifferentiability $\phi_j:=\pi/2-j\phi\ $ for $j=1..k$\\
$$a_j=h(\phi_j)$$
$$b_j+c_j=h'(\phi_j^+)-h'(\phi_j^-)=\frac{1}{2}(n-1)\frac{r^{-j}}{1-r^{-k}}(\cos'(\phi_j^+)-\cos'(\phi_j^-))=
(n-1)\frac{r^{-j}}{1-r^{-k}}$$
$$b_j-c_j=h'(\phi_j^+)+h'(\phi_j^-)=\frac{n-1}{1-r^{-k}}\sum_{i=1..k,i\neq
j}r^{-i}\frac{d}{d\alpha}|\cos(\alpha+i\phi)|_{\alpha=\phi_j}$$

\be \label{tri1}L=x_0+\bigcup_{j=1..k}e^{i\phi_j}(T(a_j,b_j,c_j)\cup
-T(a_j,b_j,c_j))\ee

Notice that if we wouldn't gather elements with the same value of
cosinus in (\ref{rat}) we can look on the construction above as(now
$j=1..\infty$):
$$a_j=h(\phi_j)$$
$$b_j+c_j=h'(\phi_j^+)-h'(\phi_j^-)=\frac{1}{2}(n-1)r^{-j}(\cos'(\phi_j^+)-\cos'(\phi_j^-))=
(n-1)r^{-j}$$
$$b_j-c_j=h'(\phi_j^+)+h'(\phi_j^-)=(n-1)\sum_{i>0,i\neq
j}r^{-i}\frac{d}{d\alpha}|\cos(\alpha+i\phi)|_{\alpha=\phi_j}$$ \be
\label{tri2} L=x_0+\bigcup_{j>0}e^{i\phi_j}(T(a_j,b_j,c_j)\cup
-T(a_j,b_j,c_j))\ee In form (\ref{tri2}) triangles from (\ref{tri1})
are constructed from infinity many triangles with disjoint interiors.\\
We will use this form to construct $L$ in $\phi\notin\pi\mathbb{Q}$ case.\\
Take $\phi^i\in\pi\mathbb{Q},\ \lim_{i\to\infty}\phi^i=\phi$ series.
We construct polygon for each $\phi^i$. It's easy to check:
$$\phi_j^i\to\phi_j\quad\quad a_j^i\to a_j\quad\quad b_j^i\to
b_j\quad\quad c_j^i\to c_j$$ for any $j>0$.\\
So (\ref{tri2}) gives convex hull in this case.\\
Its boundary contains countable number of segments.\\

We can now find a formula for length of the boundary($B$) of $L$ and
it's area($A$):
$$ B=2(n-1)\sum_{j>0}r^{-j}=2\frac{n-1}{r-1}$$ \beas
A&=&\frac{1}{2}(n-1)^2\sum_{i>0}r^{-i}\sum_{j>0}r^{-j}|\cos((i-j)\phi+\pi/2)|=
\frac{1}{2}(n-1)^2|\sin(0)|(r^{-2}+r^{-4}+...)+\\
&+&(n-1)^2\sum_{v:=|i-j|=1..\infty}|\sin(v\phi)|(r^{-2-v}+r^{-4-v}+...)=
\frac{(n-1)^2}{r^2-1}\sum_{v>0}|\sin(v\phi)|r^{-v} \eeas It's
interesting that $B$ doesn't depends on $\phi$.\\

We can get an interesting trigonometric inequality from this
formulas, namely we know (eg \cite{cia}) that for given length of
boundary, the largest area has circle, so $A\leq B^2/4\pi$:
\begin{sps}$\sum_{j>0}|\sin(j\phi)|r^{-j}\leq\frac{1}{\pi}\frac{r+1}{r-1}\quad$
for any $r>1, \phi\in[0,2\pi]$.
\end{sps}

Digression: to analyze points of indifferentiability in higher
dimension - they will correspond to pyramids, which base can be
analyzed by cutting space with two-dimensional planes - jumps of
derivatives in different direction gives the width function for the
base of pyramid. As above, we can also do it in different way - by
analyzing differentiable points around indifferentiable one.
\section{Conclusion}
We have shown that self similarity equation can be written in terms
of function describing some of its properties, like width in any
direction. Other property which can be written in that way and can
give some interesting results can be $f(x):=\mu(K\cap(K+x))$.\\
Obtained functional equation can be usually approximated numerically
and used to approximate, analyze our set.\\
The width function can be used for example to exclude some points or
to find convex hull of it and calculate some of its properties.

\end{document}